\documentclass{article}

\newtheorem{theorem}{Theorem}
\newtheorem{corollary}{Corollary}

\newtheorem{proposition}{Proposition}

\begin{document}

\title{On the homology theory of fiber spaces}

\author{T.V. Kadeishvili}

\date { }

\maketitle

\vspace{10mm}

{\bf This paper was published (in Russian) in Uspekhi Mat. Nauk
35:3 (1980), 183-188. The English translation was published in
Russian Math. Surveys, 35:3 (1980), 231-238.}

\vspace{10mm}

In this paper  the homology theory of fibre spaces is studied by introducing additional
algebraic structure in homology and cohomology.

All modules are assumed to be over an arbitrary associative ring
$\Lambda$ with unit; by a differential algebra, coalgebra, module,
or comodule we mean these objects graded by non-negative integers;
$\hat{a}$ denotes $(-1)^{deg a}$.

{\bf The category $A(\infty)$.} An {\it $A(\infty)$-algebra in the
sense of Stasheff} \cite{Sta} is defined to be a graded
$\Lambda$-module $M$, endowed with a set of operations
$\{m_i:\otimes^iM\to M, i=1,2,. . . \}$ satisfying the conditions
$m_i((\otimes^iM)q)\subset M_{q+i-2}$ and
$$
\sum_{k=0}^{i-1} \sum_{j=1}^{i-k} (-1)^k
m_{i-j+1}(\hat{a}_1\otimes . . . \otimes \hat{a}_k\otimes
m_j(a_{k+1}\otimes . . . \otimes a_{k+j})\otimes a_{k+j+1}\otimes
. . . \otimes a_i)=0
$$
for any $a_i\in M$ and $i\geq 1$. A {\it morphism of
$A(\infty)$-algebras} $ (M,\{m_i\})\to (M' ,\{m'_i\}) $ is a set
of homomorphisms $\{f_i:\otimes^iM\to M', i=1,2,. . . \}$
satisfying the conditions $f_i((\otimes^iM)q)\subset M'_{q+i-1}$
and
$$
\begin{array}{l}
\sum_{k=0}^{i-1} \sum_{j=1}^{i-k} (-1)^k
f_{i-j+1}(\hat{a}_1\otimes . . . \otimes \hat{a}_k \otimes
m_j(a_{k+1}\otimes
. . . \otimes a_{k+j})\otimes . . .  \otimes a_i)=\\
\sum_{t=1}^{i} \sum_{S(t,i)} m'_t(f_{k_1}(a_1\otimes. . . \otimes
a_{k_1})\otimes. . . \otimes f_{k_t}(a_{i-k_t+1}\otimes. . .
\otimes a_{i}))
\end{array}
$$
where $S(t,i)=\{k_1,. . . ,k_t\in N,\ \sum k_p=i\}$. The
$A(\infty)$-algebras together with these morphisms form a
category, which we denote by {\bf $A(\infty )$}.

The specification on $M$ of an $A(\infty )$-algebra structure
$(M,\{m_i\})$ is equivalent to the specification on the tensor
coalgebra $T^c(M)=\Lambda + M + M\otimes M+ . . . $ with the
grading $dim(a_1\otimes. . . \otimes a_n)=\sum dim a_i+n$ and
comultiplication
$$
\Delta (a_1\otimes. . . \otimes a_n)=\sum_{i=0}^{n}(a_1\otimes. . . \otimes
a_i)\otimes (a_{i+1}\otimes. . . \otimes a_n)
$$
of a differential $d_m:T^c(M)\to T^c(M)$ that turns $T^c(M)$ into
a differential coalgebra; this set $\{m_i\}$ determines the
differential $d_m$ by
$$
d_m(a_1\otimes. . . \otimes a_n)= \sum_{k=0}^{n-1}
\sum_{j=1}^{n-k} (-1)^k  \hat{a}_1\otimes . . . \otimes \hat{a}_k
\otimes m_j(a_{k+1}\otimes . . . \otimes a_{k+j})\otimes . . .
\otimes a_n,
$$
and the differential coalgebra $(T^c(M),d_m)$ is called the
$\tilde{B}$-construction of the $A(\infty )$-algera $(M,\{m_i\})$
(Stasheff \cite{Sta}) and is denoted by $\tilde{B}(M,\{m_i\})$.
The specification of an $A(\infty )$-algebra morphism
$\{f_i\}:(M,\{m_i\})\to (M',\{m'_i\})$   is equivalent to  that of
a differential coalgebra mapping $f:\tilde{B}(M,\{m_i\})\to
\tilde{B}(M',\{m'_i\})$; the morphism $\{f_i\}$ determines the
mapping $f$ by
$$
f(a_1\otimes. . . \otimes a_n)= \sum_{t=1}^{n} \sum_{S(t,n)}
f_{k_1}(a_1\otimes. . . \otimes a_{k_1})\otimes. . . \otimes
f_{k_t}(a_{n-k_t+1}\otimes. . . \otimes a_{i}).
$$
Thus the category {\bf $A(\infty )$} can be identified  with a
full subcategory of the category of differential coalgebras.

An arbitrary object  in {\bf $A(\infty )$} of the form
$(M,\{m_1,m_2,0,0,. . . \})$  is identified with the differential
algebra $(M,\partial,\cdot )$ where $\partial=m_1$  and $a_1\cdot
a_2=-m_2(\tilde{a}_1\otimes a_2)$. For such an object the
$\tilde{B}$-construction coincides with the usual
$B$-construction, any morphism of such objects of the form
$\{f_1,0,0,. . . \}$ is identified with the differential algebra
mapping $f_1:(M,\partial,\cdot )\to (M',\partial',\cdot)$. Thus
the category of differential algebras is a subcategory of {\bf
$A(\infty )$}, while the category {\bf $DASH$} (see \cite{Munk})
is the full subcategory of {\bf $A(\infty )$} generated by
differential algebras, and the functor $\tilde{B}$ is an extension
of  $B$ from this subcategory to {\bf $A(\infty )$}.

\begin{theorem}
\label{Thm1} For  any differential algebra $C$  with free $H_i(C
), i\geq 0$ it is possible to introduce on $H(C )$ an $A(\infty
)$-algebra structure
$$
(H(C ),\{X_i\}),\ X_i:\otimes^iH(C )\to H(C ), \ i=1,2,3,. . .
$$
such that $X_1=0,\ X_2(a_1\otimes a_2)=-\tilde{a}_1\cdot a_2$  and
there exists an $A(\infty )$-morphism
$$
\{f_i\}:(H(C ),\{X_i\})\to  (C,\{m_1,m_2,0,0,. . . \})
$$
for which $f_1:H(C )\to C$ induces an identical isomorphis in
homology.
\end{theorem}

\noindent {\bf Proof.} We need to construct two sets of homomorphisms
$$
\{X_i:\otimes^iH(C )\to H(C ),\ i=1,2,3,. . . \}, \
\{f_i:\otimes^iH(C )\to C),\ i=1,2,3,. . . \},
$$
satisfying the conditions in the definition of the category $A(\infty )$:
\begin{equation}
\label{X} \sum_{k=0}^{i-1} \sum_{j=1}^{i-k} (-1)^k
X_{i-j+1}(\hat{a}_1\otimes . . . \otimes \hat{a}_k\otimes
X_j(a_{k+1}\otimes . . . \otimes a_{k+j})\otimes a_{k+j+1}\otimes
. . .  \otimes a_i)=0,
\end{equation}

\begin{equation}
\begin{array}{l}
\label{f} \sum_{k=0}^{i-1} \sum_{j=1}^{i-k} (-1)^k
f_{i-j+1}(\hat{a}_1\otimes . . . \otimes \hat{a}_k\otimes
X_j(a_{k+1}\otimes
. . . \otimes a_{k+j})\otimes . . . \otimes  a_i)=\\
m_1f_{i}(a_1\otimes. . . \otimes a_{i})+\\
\sum_{s=1}^{i-1}m_2(f_{s}(a_1\otimes. . . \otimes a_{s})\otimes
f_{i-s}(a_{s+1}\otimes. . . \otimes a_{i}))
\end{array}
\end{equation}
for arbitrary $a_k\in H(C )$  and $i\geq 1$. Fir $i=1$ we take
$X_1=0$, and,  using  the fact that $H_i(C )$ is free, we define
$f_1:H(C )\to C$ to be a cycle-choosing homomorphis; the
conditions (\ref{X})  and (\ref{f}), as well as  the initial
condition on $f_1$ are thereby satisfied. Suppose now that $X_i$
and $f_i$ have been constructed for $i<n$ in such a way that the
conditions  (\ref{X})  and (\ref{f}) hold. Let
$$
\begin{array}{l}
U_n(a_1\otimes. . . \otimes a_n)= \sum_{s=1}^{n-1}
m_2(f_{s}(a_1\otimes. . . \otimes a_{s})\otimes
f_{n-s}(a_{s+1}\otimes. . . \otimes a_{n})) +\\
\sum_{k=0}^{n-2} \sum_{j=2}^{n-1} (-1)^{k+1}
f_{n-j+1}(\hat{a}_1\otimes . . . \otimes \hat{a}_k \otimes
X_j(a_{k+1}\otimes . . . \otimes a_{k+j})\otimes . . .  \otimes
a_n)
\end{array}
$$
(here  the $X_i$ and $f_i$ already defined are involved). Then the
condition (\ref{f}) takes the form
\begin{equation}
\label{U}
m_1f_n(a_1\otimes. . . \otimes
a_n)=(f_1X_n-U_n)(a_1\otimes. . . \otimes a_n).
\end{equation}
Direct calculations show that $\partial U_n=0$, that is, $
U_n(a_1\otimes. . . \otimes a_n)$ is a cycle in $C$ for arbitrary
$a_i\in H(C )$, and we define $ X_n(a_1\otimes. . . \otimes a_n)$ to be
the class of this  cycle, that is, $X_n=\{U_n\}$. Since $f_1$  is
a cycle-choosing homomorphism, the difference $f_1X_n-U_n$ is
homological to zero. Assuming that $a_i\in H(C )$ are free
generators we define $ f_n(a_1\otimes. . . \otimes a_n)$ as an element
of $C$ boundaring this difference and extend by linearity. For the
$X_n$ and $f_n$ thus defined the condition (\ref{U}) is
automatically satisfied. The remanning condition (\ref{X}) can be
proved by a straightforward check.

We remark that the theorem is true also when an arbitrary
$A(\infty)$-algebra is taken instead of $C$, and $H(M)$ is
understood to be the homology of $M$ with respect to the
differential $m_1$.

The $A(\infty)$-algebra $(H(C ),\{X_i\})$ we call the {\it
homology $A(\infty)$-algebra of the differential algebra $C$}. As
is clear from the proof, this structure is not uniquely determined
on $H(C )$ (there is an arbitrariness in the choice of the $f_i$).
We show later that the structure of  the homology
$A(\infty)$-algebra on $H(C )$ is unique up to isomorphism in {\bf
$A(\infty )$}.

We mention that if $a_1\cdot a_2=a_2\cdot a_3=0$ for
$a_1,a_2,a_3\in H(C )$, then $X_3(a_1\otimes a_2\otimes a_3)$ is
an element of the Massey product $<a_1,a_2,a_3>$, and this fact
provides us with examples in which the operation $X_3$ is
non-trivial. The next result follows from the Theorem \ref{Thm1}.

\begin{corollary}
\label{Cor1} The mapping of differential coalgebras
$$
f: \tilde{B}(H(C ),\{X_i\})\to B(C )
$$
induces an isomorphism in homology.
\end{corollary}

{\bf The category $M(\infty)$.} An $A(\infty)$-module over an
$A(\infty)$-algebra $(M,\{m_i\})$ we define to be a graded
$\Lambda$-module $P$, endowed with a set of operations
$\{p_i:(\otimes^{i-1}M)\otimes P\to P,\ i=1,2,3,. . . \}$
satisfying the conditions $p_i(((\otimes^{i-1}M)\otimes
P)_q)\subset M_{q+i-2}$ and
$$\begin{array}{l}
\sum_{k=0}^{i-2} \sum_{j=1}^{i-k-1} (-1)^k
p_{i-j+1}(\hat{a}_1\otimes . . . \otimes \hat{a}_k\otimes
m_j(a_{k+1}\otimes . . . \otimes a_{k+j}) \otimes . . . \otimes\\
a_{i-1}\otimes b)+ \sum_{k=0}^{i-1}  (-1)^k
p_{k+1}(\hat{a}_1\otimes . . . \otimes \hat{a}_k\otimes
p_{i-k}(a_{k+1}\otimes . . .  \otimes a_{i-1}\otimes b))=0.
\end{array}
$$

The specification on $P$ of an $A(\infty)$-module structure over
$(M,\{m_i\})$ is equivalent to the specification on
$\tilde{B}(M,\{m_i\})\otimes P$ of a differential that turns it
into a differential comodule over $\tilde{B}(M,\{m_i\})$. The
objects of the category {\bf $M(\infty)$} are defined to be the
pairs $((M,\{m_i\}),(P,\{p_i\}))$, where $(M,\{m_i\})$ is an
$A(\infty)$-algebra, and $(P,\{p_i\})$ is an $A(\infty)$-module
over it. A morphism is defined to be a pair of sets of
homomorphisms $\{f_i\},\{g_i\}$ where $\{f_i\}:(M,\{m_i\})\to
(M',\{m'_i\})$ is a morphism of $A(\infty)$-algebras and
$$
\{g_i:(\otimes^{i-1}M)\otimes P\to P',\ i=1,2,3,. . . \}
$$
is a set satisfying the conditions $g_i(((\otimes^{i-1}M)\otimes
P)_q)\subset P'_{q+i-1}$ and
$$
\begin{array}{l}
\sum_{k=0}^{i-2} \sum_{j=1}^{i-k-1}  (-1)^k \\
g_{i-j+1}(\hat{a}_1\otimes . . . \otimes \hat{a}_k\otimes
m_j(a_{k+1}\otimes . . . \otimes a_{k+j})\otimes  . . .
\otimes a_{i-1}\otimes b)= \\
\sum_{t=1}^{i} \sum_{S(t,i)}  p'_t( f_{k_1}(a_1\otimes. . .
\otimes a_{k_1})\otimes f_{k_2}(a_{k_1+1}\otimes. . . \otimes
a_{k_1+k_2})\otimes  . . .\otimes \\
f_{k_{t-1}}(a_{k_1+. . . +k_{t-2}+1}\otimes. . . \otimes a_{k_1+.
. . +k_{t-1}}) \otimes g_{k_t}(a_{i-k_t+1}\otimes. . . \otimes
a_{i-1}\otimes b));
\end{array}
$$
These conditions ensure that the mapping
$$
g:\tilde{B}(M,\{m_i\})\otimes P\to \tilde{B}(M',\{m'_i\})\otimes
P'
$$
given by
$$
\begin{array}{l}
g(a_1\otimes . . . \otimes a_{i-1}\otimes b)= \sum_{t=1}^{i}
\sum_{S(t,i)}f_{k_1}(a_1\otimes. . .\otimes a_{k_1})\otimes. . .
\otimes\\
 f_{k_{t-1}}(a_{k_1+. . . +k_{t-2}+1}\otimes. . . \otimes
a_{k_1+. . . +k_{t-1}})\otimes g_{k_t}(a_{k_1+. . .
+k_{t-1}+1}\otimes. . . \otimes a_{i-1}\otimes b)
\end{array}
$$
is a differential comodule mapping compatible with $f$. With the
obvious morphisms the category of pairs $(C,D)$, where $C$ is a
differential algebra and $D$ is a differential module over it,
forms a subcategory of {\bf $M(\infty)$}.

\begin{theorem}
\label{Thm2} If $C$ is a differential algebra and $D$ is a
differential module over it such that $H_i(C )$ and $H_i(D)$ are
free, then on $H(D)$ it is possible to introduce the structure of
an $A(\infty )$-module $(H(D),\{Y_i\}),\ Y_i:(\otimes^{i-1}H(C
))\otimes H(D)\to H(D),\ p=1,2,3,. . . $ over the homology
$A(\infty)$-algebra $(H(C ),\{X_i\})$ such that $Y_1=0,\
Y_2(a\otimes b)=-\tilde{a}\cdot b$ and there exists a morphism
$(\{f_i\},\{g_i\}):((H(C ),\{X_i\}),(H(D),\{Y_i\}))\to (C,D)$ of
{\bf $M(\infty)$} for which $f_1:H(C )\to C$ and $g_1:H(D )\to D$
induce identical isomorphisms in homology.
\end{theorem}

\noindent {\bf Proof.} The sets $\{g_i\}$ and $\{Y_i\}$ are
constructed by induction on $i$ just as in the proof of Theorem
\ref{Thm1}. Using the fact that $H_i(D)$  is free, we define $g_1$
to be a cycle-choosing homomorphism, while $Y_1=0$, and  the
conditions of the category {\bf $M(\infty)$} are satisfied for
$i=1$. Let
$$
\begin{array}{l}
V_n(a_1\otimes. . . \otimes a_{n-1}\otimes b))=
\\\sum_{s=1}^{n-1} p_2(f_{s}(a_1\otimes. . . \otimes
a_{s})\otimes
g_{n-s}(a_{s+1}\otimes. . . \otimes a_{n-1}\otimes b)) +\\
\sum_{k=0}^{n-3} \sum_{j=2}^{n-1} (-1)^{k+1}
g_{i-j+1}(\hat{a}_1\otimes . . . \otimes \hat{a}_k \otimes
X_j(a_{k+1}\otimes . . . \otimes a_{k+j})\otimes . . .  \otimes\\
a_{n-1}\otimes b)+ \sum_{k=1}^{n-2}(-1)^k g_{k+1}(\hat{a}_1\otimes
. . . \otimes \hat{a}_k \otimes Y_{n-k}(a_{k+1}\otimes . . .
\otimes a_{n-1}\otimes b)),
\end{array}
$$
then $\partial V_n=0$, therefore, we define $Y_n=\{V_n\}$. Since
$g_1$ is a cycle-choosing homomorphism, $ g_1Y_n-V_n $ is a cycle
in $D$ homologous to zero. Using the fact that $H_i(D)$ is free,
we define $g_n:(\otimes^{n-1})\otimes H(D)\to H(D)$ to be a
homomorphism for which $\partial g_n= g_1V_n-V_n$. The conditions
of the category $M(\infty)$ are satisfied for the $Y_n$ and $g_n$
thus defined.

{\bf Twisted tensor products.} The twisted tensor products of
Brown \cite{Brown} can be generalized from the case of
differential algebras and modules to the case of
$A(\infty)$-algebras and $A(\infty)$-modules: for an arbitrary
differential coalgebra $(K,d)$ and an $A(\infty)$-algebra
$(M,\{m_i\})$ a $\sim$-twisting cochain we define to be a
homomorphism $\phi:K\to M$ of degree $-1$ that satisfies the
condition
$$
\phi d=\sum_{i=1}^{\infty}m_i(\phi\otimes\phi\otimes . . . \otimes
\phi) \Delta^i ,
$$
where $\Delta^i:K\to K\otimes K$ is the homomorphism defined by $
\Delta^1=id_K,\ \Delta^2=\Delta :K\to K\otimes K,\
\Delta^i=(id_K\otimes \Delta^{i-1})\Delta$. The specification of a
$\sim$-twisting cochain $\phi:K\to M$ is equivalent to that of a
mapping of differential coalgebra $f_{\phi}: (K,d)\to
\tilde{B}(M,\{m_i\})$. For any $(K,d)$ and
$((M,\{m_i\}),(P,\{p_i\}))\in M(\infty)$ any $\sim$-twisting
cochain $\phi:K\to M$ on the tensor product $K\otimes P$
determines by
$$
\partial_{\phi}=d\otimes id_P+\sum_{i=1}^{\infty}(\hat{id}\otimes p_i)
(id_K\otimes\phi\otimes . . .  \otimes \phi\otimes id_P)(\Delta^i\otimes id_P)
$$
a differential, turning $(K\otimes P,\partial_{\phi})$ into a
differential comodule over $(K,d)$; this differential comodule is
called the {\it $\sim$-twisted tensor product} $K\otimes_{\phi}P$.
If $M$ is an $A(\infty)$-algebra of the form $(M,\{m_1,m_2,0,0,. .
.\})$, and $P$  is an $A(\infty)$-module of the form
$(P,\{p_1,p_2,0,0,. . .\})$, then $\phi$ is the usual twisting
cochain, and $K\otimes_{\phi}P$ coincides with the usual twisted
tensor product $K\otimes_{\phi}P$.

We need the concept of equivalence of twisting cochains (see
\cite{Berik1}, \cite{Berik2}, \cite{Smir}). We say that
$\phi,\psi:K\to C$ are {\it equivalent} if there is a homomorphism
$c:K\to C$ of degree $0$ for which $c_0=c|C_0=0$ and
$\psi=(1+c)\star \phi$ where
$$
(1+c)\star \phi=(1+\hat{c})\cdot \phi \cdot
(1+c)^{-1}-(cd+\partial c)\cdot (1+c)^{-1};
$$
$\phi \sim \psi$ if and only if $f_{\phi},f_{\psi}:K\to B©$ are
homotopic in the sense of \cite{Munk} (coderivation homotopy):
$f_{\phi}-f_{\psi}=\partial D+D\partial$ with $(D\otimes
f_{\psi}+\hat{f}_{\phi}\otimes D)\Delta = \Delta D$.

\begin{theorem}
\label{Thm3} If $(K,d)$ is a differential coalgebra with free
$K_i$, and $\phi:K\to C$ is an arbitrary twisting cochain, then
there exists a $\sim$-twisting cochain $\phi^*:K\to H(C )$ such
that $\phi$ and $ f^*\phi^*=\sum_{i=1}^{\infty} f_i(\phi*\otimes .
. . \otimes \phi*)\Delta^i$ are equivalent.
\end{theorem}
\noindent {\bf Proof.} To construct $\phi*$ we prove the following
inductive assertion: for any $i>0$ there exists a twisting cochain
$\phi^{(i)}:K\to C$ and a homomorphisms $\phi^*_i:K_i\to H_{i-1}(C
)$ and $c_i:K_i\to C_i$ such that
$$
\begin{array}{ll}
(a)& \phi^*_{i}d=\sum_{t=2}^{i}\sum_{S(t,i)} X_t(\phi^*_{k_1}\otimes . . .
\otimes \phi^*_{k_t})\Delta^t;\\
(b) & \phi^{(i)}=(1+c_i)\star \phi^{(i-1)};\\
(c) & \phi^{(i)}_i= \sum_{t=1}^{i}\sum_{S(t,i)}
f_t(\phi^*_{k_1}\otimes . . .\otimes \phi^*_{k_t})\Delta^t.
\end{array}
$$

For $i=1$ we take $\phi^*_1=\{\phi_1\}$. Since the difference
$(\phi_1-f_1\phi^*_1)(k)$ is homologous to zero for each $k\in
K_1$ and $K_1$ is free, we obtain a homomorphism $c_1:K_1\to C_1$
for which $-\partial c_1=\phi_1-f_1\phi^*_1$. We define
$\phi^{(1)}=(1+c_1)\star \phi$, so $\phi_1^{(1)}=\phi_1+\partial
c_1=f_1\phi^*_1$. Suppose now that $\phi^{(i)},\ \phi^*_i$, and
$c_i$ have already been constructed in such a way that $(a)$,
$(b)$, and $(c)$ hold for $i<n$. Let
$$
W_n=\phi^{(n-1)}_n -
\sum_{t=2}^{n}\sum_{S(t,n)} f_t(\phi^*_{k_1}\otimes . . .\otimes \phi^*_{k_t})\Delta^t;
$$
A direct check shows that $\partial W_n=0$; we define
$\phi^*_n=\{W_n\}$. Since the difference $W_n-f_1\phi^*_n$ is
homological to zero and $K_n$ is free, we can construct a
$c_n:K_n\to C_n$ such that $-\partial c_n= W_n-f_1\phi^*_n$; let
$\phi^{(n)}=(1+c_n)\star \phi^{(n-1)}$. Then
$$
\phi^{(n)}_n=\phi^{(n-1)}_n+\partial c_n=f_1\phi^*_n -W_n=
\sum_{t=1}^{n}\sum_{S(t,n)} f_t(\phi^*_{k_1}\otimes . . .\otimes
\phi^*_{k_t})\Delta^t,
$$
consequently, $(b)$ and $(c)$ hold for $\phi^{(n)}, \phi^*_n$, and
$c_n$, and the validity of $(a)$ can be checked directly. From
$(a)$ we see that $\phi^*=\sum_i \phi^*_i$ is a $\sim$-twisting
cochain, and from $(b)$ and $(c)$ we deduce that
$f^*\phi^*=\phi^{\infty}$, where $\phi^{\infty}= \Pi_i
(1+c_i)\star \phi \sim\phi$.

It follows from Theorem \ref{Thm3} that for any differential
coalgebra  mapping $g:K\to B©$ there exists a $G^*:K\to \tilde{B}
(H(C ),\{X_i\})$ for which $g$ and $fg^*$ are homotopic in the
sense of \cite{Munk}. This assertion implies uniqueness mentioned
above for the structure of homology $A(\infty)$-algebra: if $(H(C
),\{X_i\})$ and $(H(C ),\{X'_i\})$ are two structures of homology
$A(\infty)$-algebra on $H(C )$, then by taking $K=\tilde{B} (H(C
),\{X'_i\})$ and $g=f': \tilde{B} (H(C ),\{X'_i\})\to B(C )$, we
obtain a
$$
g^*:\tilde{B} (H(C ),\{X'_i\})\to \tilde{B} (H(C
),\{X_i\})
$$
for which $fg^*\sim g$. Then the first component of the
$A(\infty)$algebra morphism $\{g^*_i\}: (H(C ),\{X'_i\})\to (H(C
),\{X_i\})$ induced by $g^*$ is $g^*_1=id_{H(C )}$, and this
implies that $\{g^*_i\}$ is an isomorphism in {\bf $A(\infty)$}.

The next result follows from Theorems \ref{Thm1}, \ref{Thm2} and
\ref{Thm3}.

\begin{corollary}
\label{Cor2} $K\otimes_{\phi}D$ and $K\otimes_{\phi^*}H(D)$ have
isomorphic homology under the conditions of Theorems \ref{Thm1},
\ref{Thm2} and \ref{Thm3}.
\end{corollary}

The results obtained have the following applications.

The first proposition is obtained from Corollary \ref{Cor1} by
taking $C=\bar{C}_*(G)$, where $G$ is a connected topological
group such that the $H_i(G)$ are free, bearing in mind that $H(B(C
))=H_*(B_G)$.

\begin{proposition}
\label{Prop1}
The homology of the $\tilde{B}$-construction
$\tilde{B}(\bar{H}_*(G),\{X_i\})$ are isomorphic to that of
classifying space $B_G$.
\end{proposition}


The next proposition is obtained from Corollary \ref{Cor1} by
taking $C=C^*(B,b_0)$, where $B$ is a simply connected space with
free groups $H^i(B,b_0)$, and bearing in mind that $H(B(C
))=H^*(\Omega C)$.

\begin{proposition}
\label{Prop2} The homology of the $\tilde{B}$-construction
$\tilde{B}(\bar{H}^*(B,b_0),\{X_i\})$ are isomorphic to the
cohomology of the loop space $\Omega B$.
\end{proposition}

Let $\xi=(X,p,B,G)$ be a principal $G$-fibration with paracompact
base and connected $G$, let $F$ be a $G$-space, and
$\xi[F]=(E,p,B,F,G)$ the associated fiber bundle, with the
$H_i(G)$ and $H(_i(F)$ free. The final proposition is obtained
from Corollary \ref{Cor2} ty taking $C=C_*(G),\ D=C_*(F)$, and
$\phi$ a twisting cochain of the fibration $\xi$ (\cite{Brown}).

\begin{proposition}
\label{Prop3} The homology of the $\sim$-twisted tensor product
$C_*(B)\otimes_{\phi^*}H_*(F)$ is isomorphic to that of $E$.
\end{proposition}

This proposition generalize a result of Shih \cite{Shih}: if $G$
is $(n-1)$-connected, then the components $\phi^*_i\in
C^i(B,H_{i-1}(G))$ vanish for $0<i<n+1$, therefore, the
differentials $d^i$ of the spectral sequence of $\xi [F]$ are
trivial for $1<i<n+1$, and the components are cocycles for
$n<i<2n+1$, consequently, $d^j$ can be expressed for $n<i<2n+1$ in
terms of certain characteristic classes of $\xi$ and the operation
$Y_2:H_*(G)\otimes H_*(F)\to H_*(F)$; we remark that higher
operations $Y_i$ are needed for computing $d^j,\ j>2n$, in terms
of $\phi^*$.

Theorems \ref{Thm1}, \ref{Thm2} and Proposition \ref{Prop3} were
announced in \cite{Kade1} and Proposition \ref{Prop2} in
\cite{Kade2}.


\begin{thebibliography}{99}



\bibitem{Sta} J.D. Stasheff, Homotopy associativity of H-spaces, I, II,
Trans. Amer. Math. Soc. 108 (1963), 27-313. MR 28  1623.


\bibitem{Munk} H. Munkholm, The Eilenberg-More spectral sequence and strongly
homotopy multiplicative maps,
J. Pure and Appl. Algebra 5 (1974), 1-50.
MR 50  3227.

\bibitem{Brown} E. H. Brown, Twisted tensor products, I, Ann. Of Math.
(2) 69 (1959), 223-246. MR 21  4423.

\bibitem{Berik1} N. A. Berikashvili, The differentials of a spectral sequence,
Bull. Of Georg. Acad. Sci., 51 (1968), 9-14. MR 41  9258.

\bibitem{Berik2} N. A. Berikashvili, Homology theory of spaces,
Bull. Of Georg. Acad. Sci., 86 (1977), 529-532. MR 57  13949.


\bibitem{Smir} V. A. Smirnov, The functor $D$ for twisted tensor products,
Mat. Zametki 20 (1976), 465-472. MR 55  4172.

\bibitem{Shih} W. Shih, Homologie des espaces fibres, Inst. Hautes Etudes Sci.
Publ. Math. 1962, no. 13, 88. MR 26  1893.

\bibitem{Kade1} T. Kadeishvili, On the differentials of spectral sequence of a fiber bundle,
Bull. Of Georg. Acad. Sci.,
82 (1976), 285-288. MR 55  6430.

\bibitem{Kade2} T. Kadeishvili, On the homology of classifying spaces,
Proc. 7th All Union Topology Conference, Minsk 1977.

\end{thebibliography}
\end{document}